\documentclass[12pt]{article}
\usepackage{graphicx}
\usepackage{ucs}
\usepackage{amsmath}
\usepackage{sectsty}
\usepackage{amssymb}
\usepackage{vmargin}
\usepackage{amsthm}
\usepackage{mathtext}
\usepackage{eufrak}
\usepackage{yfonts}
\usepackage{amsthm}

\usepackage[mathscr]{eucal}
\usepackage{amsbsy}
\usepackage{amsmath}
\usepackage{graphics}
\usepackage{epsfig}

\usepackage{latexsym}
\setpapersize{A4}
\setmarginsrb{2.5cm}{2cm}{1.5cm}{2cm}{0pt}{0mm}{0pt}{13mm}
\usepackage{indentfirst}
\newtheorem{lemma}{Lemma}
\newtheorem{remark}{Remark}
\newtheorem{theorem}{Theorem}
\newtheorem{corollary}{Corollary}
\newcommand{\Var}{\mathrm{Var}}
\newcommand{\E}{\mathrm{E}}

\newcommand{\halmos}{\vspace{3mm} \hfill \mbox{$\Box$}\\[2mm]}

\date{}

\begin{document}

\title{Random walks in non homogeneous Poissonian environment}
\vspace{15pt}

\author{Youri Davydov $^{\text{\small 1}}$  and  Valentin Konakov$^{\text{\small 2,}}$ 
\footnote{For the second author the article was prepared within the framework of a subsidy granted to the HSE by the Government of the Russian Federation for the implementation of the Global Competitiveness Program.}\\
{\small $^{\text{1}}$ Universit\'e Lille 1
%, Laboratoire Paul Painlev\'e, \; \;
and Saint Peterbourg State University}\\
 {\small $^{\text{2}}$  Higher Scool of Economics}\\}
 %{\small and Institute of Mathematics and Informatics,  }}
%\institute{University of Lille 1, France }

\vspace{20pt}

\maketitle

\section{Introduction}

 We consider the moving particle process in $R^{d}$ which is defined in
the following way. There are two independent sequences 
$\left( T_{k}\right) $
and $\left( \varepsilon _{k}\right) $ of random variables. 

   The variables $T_k$ are non negative 
  and $\forall k \;\;T_k \leq T_{k+1},$ while variables 
  $\varepsilon _{k}$ form an i.i.d sequence with
common distribution concentrated on the unit sphere $S^{d-1}.$

    The values 
$\varepsilon _{k}$ are interpreted as the directions, and $T_{k}$ as the moments
of change of directions.

 A particle starts from zero and moves in the
direction $\varepsilon _{1}$ up to the moment $T_{1}.$ It then changes direction
to $\varepsilon _{2}$ and moves on within the time interval 
$T_{2}-T_{1},$ etc.
The speed is constant at all sites. The position of the particle at time $t$
is denoted by $X(t).$
    
    Study of the processes of this type has a long history.
 The first work dates back probably to Pearson  and continued by Kluyer (1906) and Rayleigh (1919).
  Mandelbrot (1982) 
considered the case where the increments $T_{n}-T_{n-1}$ form i.i.d.
sequence with the common law having a heavy tail. He also introduced the
term "Levy flights" later changed to  "Random flights". 

  To date, a
large number of works were accumulated, devoted to the study of such
processes, we mention here only articles by Kolesnik (2009),
 Orsingher and  De Gregorio
(2012, 2015) and  Orsingher and  Garra (2014) which contain an extensive bibliography and where for different
assumptions on $\left( T_{k}\right) $ and $\left( \varepsilon _{k}\right) $ the
exact formulas for the distribution of $X(t)$ were derived. 

Our goals are different. 
 
   Firstly, we are interested in the global behavior
of the process $X=\{X(t),\;t\in R_{+}\},$ namely, we are looking for conditions
under which the  processes $\{Y_{T},\;T>0\},$

\[
Y_{T}(t)=\frac{1}{B(T)}X(tT),\;\;t\in \lbrack 0,1],
\]%
weakly converges in $C[0,1]:$ $\;Y_{T}\Longrightarrow Y,\; B_T\longrightarrow \infty,\; T\longrightarrow
\infty .$

 From now on we suppose that the points $(T_k),\; T_k\leq T_{k+1},$
    form a Poisson point process  in $R_{+}$ denoted by 
    $\mathbf{T}.$ 
\vspace{5pt}
    
 It is clear that in the homogeneous case the process $X(t)$ is a conventional
random walk because the spacings $T_{k+1}-T_k$ are  independent, and then the limit process is Brownian motion. 
\vspace{5pt}

In the non homogeneous case the situation is more complicated as
 these spacings  are not independent. Nevertheless
it was possible to distinguish three modes that
determine different types of limiting processes. 
\vspace{5pt}

 For a more precise description of
the results it is convenient to assume that $T_{k}=f(\Gamma _{k}),$ where 
$ \mathbf{\Pi}= \left( \Gamma _{k}\right) $ is a standard homogeneous Poisson point process on $R_{+}$ with intensity $1.$ 
In this case
\[
(\Gamma_k)\stackrel{\mathcal{L}}{=}(\gamma_1+\gamma_2+\dots+\gamma_k)
\]
where $(\gamma_k)$ are i.i.d standard exponential random variables.

If the function $f$  has power growth, 
$$
f(t)=t^{\alpha },\,\alpha \geq 1,
$$ 
the behavior of the
process is analogous to the uniform case and then in the limit
we obtain a Gaussian process which is a linearly transformed Brownian motion
%Gaussian process%
\[
Y(t)=\int_{0}^{t}K_\alpha(s) dW(s),
\]%
where $W$ \ is a  process of Brownian motion, for which the covariance matrix of $W(1)$
coincides with  the covariance matrix of $\varepsilon _{1}$
and $K_\alpha(s)$ is a nonrandom kernel, its exact expression is given below. 

 In the case of exponential growth, 
$$
f(t)=e^{t\beta },\;\beta >0,
$$
 the
limiting process is piecewise linear with an infinite number of units, but 
$\forall \epsilon >0$ the number of units in the interval 
$[\varepsilon, 1]$
will be  a.s. finite.
%\vspace{5pt}

Finally, with the super exponential growth of $f$, the process degenerates: its
trajectories are linear functions:

\[
Y(t)=\varepsilon t,\;\text{ }t\in \lbrack 0,1],\;\text{ }\varepsilon \overset{Law}{=}%
\varepsilon _{1}.
\]%
\vspace{7pt}

%%%%%%%%%%%%%%%%%%%%%%%%%%%%%%%%%%%%%%%%%%%%%%%%%
   In the second part of the paper the process $X(t)$ is considered as a Markov
chain. We construct diffusion approximations for this process and
investigate their accuracy. To prove the  weak convergence we use the
approach of Stroock and Varadhan (1979). Under our assumptions the diffusion
coefficients $a$ and $b$ have the property that for each $x\in R^{d}$ the
martingale problem for $a$ and $b$ has exactly one solution $P_{x}$ starting
from $x$ (that is well posed). It remains to check the conditions from Stroock and Varadhan (1979)
which imply the weak convergence of our sequence of Markov chains to this
unique solution $P_{x}$ . We consider also more general model which may be
called as "random walk over ellipsoids in $R^{d}$ ". For this model we
establish the convergence of transition densities and obtain Edgeworth type
expansion up to the order $n^{-3/2}$, where $n$ is a number of switching.
The main tool in this part is the paramertix method (Konakov (2012), Konakov and Mammen (2009)).      

%%%%%%%%%%%%%%%%%%%%%%%%%%%%%%%%%%%%%%%%%%%%%%%%

\section{Random flights in Poissonian environment}

The reader is reminded that we suppose $T_{k}=f(\Gamma _{k}),$ where $%
\left( \Gamma _{k}\right) $ is a standard homogeneous Poisson point process
on $R_{+}.$   Assume also that $E\varepsilon_1 = 0.$    

 It is more convenient  to consider at first the behavior of processes 
$$
Z_n(t) = Y_{T_n}(t),
$$
as for $T=T_n$ the paths of $Z_n$  have an 
integer number of full  segments on the interval [0,1]. 
The typical path of $\{Z_n(t), t\in[0,1]\}$  is a continuous broken line  with vertices\\
 $\{(t_{n,k},\;\frac{S_k}{B_n}),\;\; k=0,1,\dots,n\}$, 
 where $\;t_{n,k}=\frac{T_k}{T_n},\; T_{0}=0,\; B_n = B(T_n),\; S_k =\sum_1^k\varepsilon_i(T_i-T_{i-1})$.
    
\begin{theorem}
 \label{th1} 
 Under previous assumptions
\begin{itemize}
\item[1)]  If 
the function $f$  has power growth: 
$\;f(t)=t^{\alpha },\,\alpha > 1/2,$ 
we take $B(T)= T^{\frac{2\alpha-1}{2\alpha}}.$

Then $Z_{n}\Longrightarrow Y,$ where $Y$ is
a Gaussian process%
\[
Y(t)=\sqrt{2\alpha }\int_{0}^{t}s^{\frac{\alpha -1}{2\alpha }}dW(s),
\]%
and $W$  is a  process of Brownian motion, for which the covariance matrix of $W(1)$
coincides with  the covariance matrix of $\varepsilon _{1}.$

\item[2)] 
If the function $f$  has exponential growth:
$\;f(t)=e^{t\beta },\;\beta >0,$
we take $B(T)= T.$
\vspace{7pt}

Then $Z_{n}\Longrightarrow Y,$ where $Y$ is a continuous piecewise
linear process with the vertices at the points $(t_k,Y(t_k)),$
$$
t_k = e^{-\beta\Gamma_{k-1}}, \;\; \Gamma_0 =0,
$$
$$
Y(t_k) = \sum_{i=k}^\infty 
\varepsilon_k(e^{-\beta\Gamma_{i-1}}-e^{-\beta\Gamma_{i}}),\;\;\;
Y(0)=0.
$$
      
\item[3)] In super exponential case,  suppose that $f$ is increasing absolutely continuous and such that 
$$
\lim_{t\rightarrow \infty} \frac{f^\prime(t)}{f(t)} = +\infty.
$$
We take $B(T)= T.$
\vspace{7pt}

Then 
$\frac{T_n}{T_{n+1}}\rightarrow 0$ in probability,  and
$Z_{n}\Longrightarrow Y,$  where the limiting process $Y$
 degenerates:
$$
Y(t)=\varepsilon_1 t,\;\text{ }t\in \lbrack 0,1]
.   
$$
$$
\,
$$
\end{itemize}

\end{theorem}

%{ Comments.}
\begin{remark}
In the case of power growth the limiting process admits  the following representation:
$$
Y(t) \stackrel{\mathcal{L}}{=} \alpha\sqrt{\frac{2}{2\alpha-1}}W(t^{\frac{2\alpha -1}{\alpha }}),
$$
where, as before, $W$ is a Brownian motion, for which the covariance matrix of $W(1)$
coincides with  the covariance matrix of $\varepsilon _{1}.$

It is clear that we can also express $Y$  in another way :
$$
Y(t) \stackrel{\mathcal{L}}{=} \alpha\sqrt{\frac{2}{2\alpha-1}}K^{\frac{1}{2}}w(t^{\frac{2\alpha -1}{\alpha }}),
$$
where  $w$  is a standard Browniam motion and $K$ is the covariance matrix of $\varepsilon_1.$
\end{remark}

\begin{remark}
In the case of exponential growth it is possible to describe the limiting process $Y$ in the following way:

We take a p.p.p. ${\mathbf T}=(t_k)$, $t_k=e^{-\beta\Gamma_{k-1}}$, defined on $(0,1]$, and define a step process\\
 $\{Z(t),\, t\in(0,1]\}$,
\[
Z(t)=\varepsilon_k\quad\text{for }\;t\in(t_{k+1},t_k].
\]

Then
\[
Y(t)=\int_0^tZ(s)\,\mathrm{d}s.
\]
\end{remark}

\section{Diffusion approximation} 
In this section firstly we consider a model of random flight
which is equivalent to the study of random broken lines
$\{X_n(t), \;t \in [0,1]\}$ with the vertices 
$(\frac{k}{n},\;X_n(\frac{k}{n})),$
and such that ($h = \frac{1}{n}$)
\[
X_n\left( (k+1)h\right) =X_n(kh)+hb(X_n(kh))+\sqrt{h}\xi_k (X(kh)),
\]%
\begin{equation}
X_n(0)= x_0,\;\;\;
\xi_k (X_n(kh))=\rho _{k}\sigma (X_n(kh))\varepsilon _{k},  \label{13}
\end{equation}%
where $\left\{ \varepsilon _{k}\right\} $ and $\left\{ \rho _{k}\right\} $ are two independent sequences and

$\left\{ \varepsilon _{k}\right\} $ are i.i.d. r. v. uniformly distributed on the unit sphere $S^{d-1},$ 

 $\left\{ \rho _{k}\right\} $ are i.i.d. r. v. having a density, $\rho _{k}\geq 0,$ $E\rho_{k}^{2}=d,$
 
$b:R^{d}\longrightarrow R^{d}$ \ 
is a bounded measurable function and  
$\sigma :R^{d}\longrightarrow R^{d}\times R^{d}$ 
is a bounded measurable matrix function.  

 \begin{theorem}\label{th2}
Let $X = \{X(t), \;t\in [0,1]\}$ be a solution of stochastic equation
$$
X(t) = x_0 +\int_0^t b(X(s))ds + \int_0^t \sigma(X(s))dw(s).
$$
Suppose that $b$ and $\sigma$ are continuous functions satisfying
Lipschitz condition
$$
|b(t)-b(s)|+|\sigma(t) - \sigma(s)| \leq K|t-s|.
$$
Moreover it is supposed that $\;b(x)$ and $\frac{1}{\det{(\sigma(x))}}$
are bounded.

Then, 
\begin{center}
$X_n \Longrightarrow X\;\;\;$
 in $\;\;\;\mathbb{C}[0,1].$
 \end{center}
\end{theorem}
\vspace{10pt}
%%%%%%%%%%%%%%%%%%%%%%%%%%%%%%%%%%%%%%%%%%%%%%%%%%%%%%%%%%%%%%%

Our next result is about approximation of transition density.
We consider now more general models given by a triplet $(b(x),\sigma (x),f(r;\theta )),\;x\in
R^{d},\;r\geq 0,\;\theta \in R^{+},$ where $b(x)$ is a
vector field, $\sigma (x)$ is a $d\times d$ matrix, $a(x):=\sigma \sigma
^{T}(x)>\delta I,\;\delta >0,$ and $f(r;\theta )$ is a radial density
depending on a parameter $\theta $ controlling the frequency of changes of
directions, namely, the frequency increases when $\theta $ decreases.
Suppose $X(0)=x_{0}$. The vector $b(x_{0})$ acts shifting a particle from $%
x_{0}$ to $x_{0}+\Delta (\theta )b(x_{0}),$ where $\Delta (\theta
)=c_{d}\theta ^{2},$ $c_{d}>0.$ Several examples of such functions $\Delta
(\theta )$ for different models will be given below. Define%
\[
\mathcal{E}_{x_{0}}(r):=\{x:\left\vert a^{-1/2}(x_{0})(x-x_{0}-\Delta
(\theta )b(x_{0}))\right\vert ^{2}=r^{2}\}, 
\]%
\[
\mathcal{S}_{x_{0}}^{d}(r):=\{y:\left\vert y-x_{0}-\Delta (\theta
)b(x_{0})\right\vert ^{2}=r^{2}\}. 
\]%

%%%%%%%%%%%%%%%%%%%%%%%%%%%
The initial direction is defined by a random variable
 $\xi _{0}$, the law of $\xi _{0}$ is a pushforward of the spherical measure on $\mathcal{S}_{x_{0}}^{d}(1)$ under affine change of variables
 \[
 x-x_0-\Delta(\theta)b(x_0)=a^{1/2}(x_0)(y-x_0-\Delta(\theta)b(x_0)).
 \]
%%%%%%%%%%%%%%%%%%%%%%%%%%%%

Then particle moves along the ray $l_{x_{0}}$
corresponding to the directional unit vector 
\[
\varepsilon _{0}:=\frac{\xi _{0}-x_{0}-\Delta (\theta )b(x_{0})}{\left\vert
\xi _{0}-x_{0}-\Delta (\theta )b(x_{0})\right\vert }, 
\]%
and changes the direction in $(r,r+dr)$ with p\bigskip robability%
\begin{equation}
\det \left( a^{-1/2}(x_{0})\right) \cdot f(r\left\vert
a^{-1/2}(x_{0})e_{0}\right\vert )dr.  \label{1}
\end{equation}%
Let $\rho _{0}$ be a random variable independent on $\xi _{0}$ and
distributed on $l_{x_{0}}$ with the radial density (\ref{1}). We consider
the point $x_{1}=x_{0}+\Delta (\theta )b(x_{0})+\rho _{0}\varepsilon _{0}.$\;
Let $(\varepsilon_k,  \rho_k)$ be independent copies of $(\varepsilon_0, \rho_0).$
Starting from $x_{1}$ we repeat the previous construction to obtain $%
x_{2}=x_{1}+\Delta (\theta )b(x_{1})+\rho _{1}\varepsilon _{1}.$ After $n$
switching we get a point $x_{n},$ 
\[
x_{n}=x_{n-1}+\Delta (\theta )b(x_{n-1})+\rho _{n-1}\varepsilon _{n-1}. 
\]%
To obtain the one-step characteristic function $\Psi _{1}(t)$ we make use of
formula (6) from  Yadrenko (1980):
\[
\Psi _{1}(t)=Ee^{i\left\langle t,\rho _{0}\varepsilon _{0}\right\rangle
}=\int_{0}^{\infty }\int_{\mathcal{E}_{x_{0}}(r)}e^{i\left\langle
t,a^{1/2}(x_{0})a^{-1/2}(x_{0})\xi \right\rangle }\mu _{\mathcal{E}%
_{x_{0}}(r)}(d\xi )d\Phi _{\mathcal{E}}(r)= 
\]%
\[
=\int_{0}^{\infty }\int_{\mathcal{S}_{x_{0}}^{d}(r)}e^{i\left\langle
a^{1/2}(x_{0})t,y\right\rangle }\lambda _{r}^{d}(dy)f(r;\theta )dr= 
\]%
\begin{equation}
=2^{\frac{d-2}{2}}\Gamma \left( \frac{d}{2}\right) \int_{0}^{\infty }\frac{J_{%
\frac{d-2}{2}}(r\left\vert a^{1/2}(x_{0})t\right\vert )}{\left( r\left\vert
a^{1/2}(x_{0})t\right\vert \right) ^{\frac{d-2}{2}}}f(r;\theta )dr,
\label{2}
\end{equation}%
where $J_{\nu }(z)$ is the Bessel function and $d\Phi _{\mathcal{E}}(r)$ is
the $F$ - measure of the layer between $\mathcal{E}_{x_{0}}(r)$ and $%
\mathcal{E}_{x_{0}}(r+dr),$ $F$ is the law of $\ \rho _{0}\varepsilon _{0}.$
Now we make our main assumption about the radial density:

{\bf (A1)} \;\;The funciton $f(r;\theta )$ is homogenious of degree $-1,$
that is
\[
f(\lambda r;\lambda \theta )=\lambda ^{-1}f(r;\theta ),\;\forall \lambda \neq
0. 
\]%
Denote by $p_{\mathcal{E}}(n,x,y)$ the transition density after $n$
switching in the RF-model described above. To obtain the one step transition
density $p_{\mathcal{E}}(1,x,y)$ (we write $(x,y)$ instead of $(x_{0},x_{1})$%
) we use the inverse Fourier transform, (\ref{2}) and \textbf{(A1). }After
easy calculations we get \textbf{\ }%
\begin{equation}
p_{\mathcal{E}}(1,x,y)=\Delta ^{-d/2}(\theta )q_{x}\left(\frac{y-x-\Delta (\theta
)b(x)}{\sqrt{\Delta (\theta )}}\right),  \label{3}
\end{equation}%
where 
\begin{equation}
q_{x}(z)=\frac{2^{\frac{d-2}{2}}\Gamma \left( \frac{d}{2}\right) }{\left(
2\pi \right) ^{d}}\int_{R^{d}}\cos \left\langle \tau ,z\right\rangle \left[
\int_{0}^{\infty }\frac{J_{\frac{d-2}{2}}(\rho \left\vert a^{1/2}(x)\tau
\right\vert )}{\left( \rho \left\vert a^{1/2}(x)\tau \right\vert \right) ^{%
\frac{d-2}{2}}}f(\rho ;c_{d})d\rho \right] d\tau .  \label{4}
\end{equation}%
Consider two examples.
\vspace{5pt}

\textbf{Example 1. }We put $\Delta (\theta )=(d+1)^{2}\theta ^{2}$ and%
\[
f(r;\theta )=\frac{1}{\Gamma (d)}r^{-1}\left( \frac{r}{\theta }\right)
^{d}\exp \left( -\frac{r}{\theta }\right) . 
\]%
Using (\ref{2}), formula 6.623 (2) on page 726 from Gradshtein and
Ryzhik (1963), and the doubling formula for the Gamma function we obtain 
\[
p_{\mathcal{E}}(1,x,y)=\Delta ^{-d/2}(\theta )q_{x}\left(\frac{y-x-\Delta (\theta
)b(x)}{\sqrt{\Delta (\theta )}}\right), 
\]%
where 
\[
q_{x}(z)=\frac{(d+1)^{d/2}}{2^{d}\pi ^{(d-1)/2}\Gamma \left( \frac{d+1}{2}%
\right) \left\vert \det a^{1/2}(x)\right\vert }e^{-\sqrt{d+1}\left\vert
a^{-1/2}(x)z\right\vert }. 
\]%
It is easy to check that%
\[
\int z_{i}q_{x}(z)=0,\;\;\;\int z_{i}z_{j}q_{x}(z)dz=a_{ij}(x). 
\]
\vspace{7pt}

\textbf{Example 2. }We put $\Delta (\theta )=\theta ^{2}/2$ and%
\[
f(r;\theta )=C_{d}r^{-1}\left( \frac{r}{\theta }\right) ^{d}\exp \left( -%
\frac{r^{2}}{\theta ^{2}}\right) , 
\]%
where $C_{d}=\frac{2^{(d+1)/2}}{(d-2)!!\sqrt{\pi }}$ \;if $d$ is odd, and $%
C_{d}=\frac{2}{\left[ (d-2)/2\right] !}$ if $d$ is even. From (\ref{2}) and
formula 6.631 (4) on page 731 
of Gradshtein and Ryzhik (1963)
we obtain%
\[
p_{\mathcal{E}}(1,x,y)=\Delta ^{-d/2}(\theta )\phi _{x}(\frac{y-x-\Delta
(\theta )b(x)}{\sqrt{\Delta (\theta )}}), 
\]

where 
\[
\phi _{x}(z)=\frac{1}{\left( 2\pi \right) ^{d/2}\sqrt{\det a(x)}}\exp \left(
-\frac{1}{2}\left\langle a^{-1}(x)z,z\right\rangle \right) .
\]%
\vspace{10pt}

It is easy to see that the transition density (\ref{3}) corresponds to the
one step transition density in the following Markov chain model%
\[
X_{(k+1)\Delta \left( \theta \right) }=X_{k\Delta \left( \theta \right)
}+\Delta \left( \theta \right) b(X_{k\Delta \left( \theta \right) })+\sqrt{%
\Delta (\theta )}\xi _{(k+1)\Delta \left( \theta \right) },
\]%
where the conditional density (under $X_{k\Delta \left( \theta \right) }=x$)
of the innovations $\xi _{(k+1)\Delta \left( \theta \right) }$ is equal to $%
q_{x}(\cdot ).$\;\;If we put $\theta =\theta _{n}=\sqrt{\frac{2}{n}},$ then $%
\Delta \left( \theta _{n}\right) =\frac{1}{n}$ and we obtain a sequence of
Markov chains defined on an equidistant grid%
\begin{equation}
X_{\frac{k+1}{n}}=X_{\frac{k}{n}}+\frac{1}{n}b(X_{\frac{k}{n}})+\frac{1}{%
\sqrt{n}}\xi _{\frac{k+1}{n}},\;\;\;X_{0}=x_{0}.  \label{5}
\end{equation}%
Note that the triplet $(b(x),\sigma (x),f(r;\theta )),\;x\in R^{d},\;r\geq
0,\;\theta \in R^{+},$ of the Example 2 corresponds to the classical Euler
scheme for the $d$-dimensional SDE%
\begin{equation}
dX(t)=b(X_{t})dt+\sigma (X_{t})dW(t),\;\;\;X(0)=x_{0}.  \label{6}
\end{equation}%
Let $p(1,x,y)$ be transition density from $0$ to $1$ in the model (\ref{6}).
We make the following assumptions
\vspace{7pt}

\textbf{(A2) }The function $a(x)=\sigma \sigma ^{T}(x)$ is uniformly
elliptic.
\vspace{7pt}

\textbf{(A3) }The functions $b(x)$ and $\sigma (x)$ and their derivatives up
to the order six are continuous and bounded uniformly in $x.$ The 6-th
derivative is globally Lipschitz. 

\begin{theorem}\label{th3}
Under  assumptions  (A2), (A3) we have
the following expansion:  for any positive integer $S$ as $n\rightarrow \infty$

\begin{equation}
\sup_{x,y\in R^{d}}\left( 1+\left\vert y-x\right\vert ^{S}\right) \cdot 
\left\vert p_{\mathcal{E}}(n,x,y)-p(1,x,y)-\frac{1}{2n}p\otimes \left(
L_{\ast }^{2}-L^{2}\right) p(1,x,y)\right\vert =O(n^{-3/2}),  \label{7}
\end{equation}%
where 
\begin{equation}
L=\frac{1}{2}\sum_{i,j=1}^{d}a_{ij}(x)\partial
_{x_{i}x_{j}}^{2}+\sum_{i=1}^{d}b_{i}(x)\partial _{x_{i}}.  \label{8}
\end{equation}%
\end{theorem}

The operator  $L_{\ast }$  in (\ref{7}) is the same operator as in (\ref{8}) but with coefficients
"frozen" ay $x.$ Clearly, $L=L_{\ast }$ but, in general, $L^{2}\neq L_{\ast
}^{2}.$ The convolution type binary operation $\otimes $ is defined for
functions $f$ and $g$ in the following way%
\[
\left( f\otimes g\right)
(t,x,y)=\int_{0}^{t}ds\int_{R^{d}}f(s,x,z)g(t-s,z,y)dz.
\]

{\bf Proof.} It follows immediately from  Theorem 1 of
Konakov and Mammen (2009).

%%%%%%%%%%%%%%%%%%%%%%%%%%%%%%%%%%%%%%%%%%%%%%%%%%%%%%%%%%%%%%%

\section{Proof of Th. 1}

\subsection{Asymptotic behaviour in case 3)}
We have, taking $B_n = B(T_n) =T_n$ :
\[
\sup_{t\in\left[0,\frac{T_{n-1}}{T_n}\right]}\|X_n(t)\|_{\infty}\le\sum_{k=1}^{n-1}\frac{T_k-T_{k-1}}{T_n}=\frac{T_{n-1}}{T_n}\xrightarrow[n\rightarrow\infty]{}0 \text{  a.s.}
\]

At the same time,
\[
X_n(1)=\frac{S_{n-1}+\varepsilon_n(T_n-T_{n-1})}{T_n}=\varepsilon_n+o(1)\Rightarrow\mathcal{P}_{\varepsilon_1}
\]

Therefore the process $X_n$ converges weakly to the process $\{Y(t)\}$, $Y(t)=\varepsilon_1t$, $t\in[0,1]$.

This process is in some sense degenerate. Hence this case is not very interesting.
%%%%%%%%%%%%%%%%%%%%%%%%%%%%%%%%%%%%%%%%%%%%%%%%%%

\subsection{Asymptotic behaviour in case 2)}

%Now we consider the case 2).

Take $B_n=T_n$ and show that the limit process $Y$ is not trivial. For simplicity fix $\beta=1$. We have now 
$t_{n,k} : = \frac{T_k}{T_n}=e^{-(\Gamma_n-\Gamma_k)}=e^{-(\gamma_{k+1}+\dots+\gamma_n)}$, and
\[
X_n\left({t_{n,k}}\right)=\sum_{i=1}^k\varepsilon_i(e^{-(\gamma_{i+1}+\dots+\gamma_n)}-e^{-(\gamma_{i}+\dots+\gamma_n)}),\quad k=1,\dots,n.
\]

The process $X_n$ is completely defined by 2 independent vectors $(\varepsilon_1,\dots,\varepsilon_n)$ and $(\gamma_1,\dots,\gamma_n)$. Hence its distribution will be the same if we replace these vectors by $(\varepsilon_n,\dots,\varepsilon_1)$ and $(\gamma_n,\dots,\gamma_1)$.
In another words, the process $(X_n(\cdot))\stackrel{\mathcal{L}}{=}(Y_n(\cdot))$, where $Y_n(\cdot)$ is a broken line with vertices $(\tau_{n,k},Y_n(\tau_{n,k}))$, $(\tau_{n,k})\downarrow$, $\tau_{n,1}=1$, $\tau_{n,k}=e^{-(\gamma_1+\dots+\gamma_{k-1})}$, $k=2,\dots,n,$ and
\[
Y_n(\tau_{n,k})=\sum_{i=k}^{n-1}\varepsilon_i\left(e^{-(\gamma_1+\dots+\gamma_{i-1})}-e^{-(\gamma_1+\dots+\gamma_{i})}\right)+\varepsilon_ne^{-(\gamma_1+\dots+\gamma_{n-1})};
\]
$Y_n(0)=0$, and $\gamma_0:=0$.

Using the notation $\Gamma_k=\gamma_1+\dots+\gamma_k$ we get the more compact formula:
\[
Y_n(\tau_{n,k})=\sum_{i=k}^{n-1}\varepsilon_i\left(e^{-\Gamma_{i-1}}-e^{-\Gamma_{i}}\right)+\varepsilon_ne^{-\Gamma_{n-1}}.
\]

Consider now the process $\{Y(t), t\in[0,1]\}$ defined as follows:
\begin{equation}
Y(0)=0,\quad Y(t_k)=\sum_{i=k}^{\infty}\varepsilon_i\left(e^{-\Gamma_{i-1}}-e^{-\Gamma_i}\right);
\end{equation}
where $t_k=e^{-\Gamma_{k-1}}$, $k=2,3,\dots$, $t_1=1$; for $t\in[t_{k+1},t_k]$ $Y(t)$ is defined by linear interpolation. The paths of $Y$ are continuous broken lines, starting at 0 and having an infinite number of segments in the neighborhood of zero.

The evident estimation
\[
\begin{aligned}
&\sup_{t\in[0,1]}|Y(t)-Y_n(t)|\le\left|\sum_{i=n}^{\infty}\varepsilon_i\left(e^{-\Gamma_{i-1}}-e^{-\Gamma_{i}}\right)\right|+e^{-\Gamma_{n-1}}\le\\
&\le\sum_{i=n}^{\infty}\left(e^{-\Gamma_{i-1}}-e^{-\Gamma_{i}}\right)+e^{-\Gamma_{n-1}}=
2e^{-\Gamma_{n-1}}\longrightarrow0\quad\text{a.s.}
\end{aligned}
\]
shows that a.s. $Y_n(\cdot)\xrightarrow{\mathbb{C}[0,1]}Y(\cdot)$.

Conclusion: In case 2),  the process $X_n$ converges weakly to $Y(\cdot)$.

\begin{remark}
In the case where $\beta\neq1$ it is simply necessary  replace $e^{-\Gamma_k}$ by $e^{-\frac{\Gamma_k}{\beta}}$.
\end{remark}
\begin{remark}
It seems that the last result could be expanded by considering more general sequences $(\varepsilon_k)$.

Interpretation: $\frac{\varepsilon_k}{|\varepsilon_k|}$ defines direction, $|\varepsilon_k|$ defines the velocity of deplacement in this direction on the step $S_k$.
\end{remark}

\subsection{Asymptotic behaviour in case of power growth}
In this case $T_k=\Gamma^{\alpha}_k$, $\;\alpha > 1/2$, $\;t_{n,k}=\frac{T_k}{T_n}=\left(\frac{\Gamma_k}{\Gamma_n}\right)^{\alpha}$, and
\begin{equation}
X_n(t_{n,k})=\frac{1}{B_n}\sum_{i=1}^k\varepsilon_i(\Gamma_i^{\alpha}-\Gamma_{i-1}^{\alpha});\,\, \Gamma_0=0,\,k=0,1,\dots,n.
\end{equation}

Let $x\in {\mathbb{R}^d}$ be such that 
$|x| = 1.$ We will show below  that
\[
\Var\left(\sum_{i=1}^n
\langle\varepsilon_i,x\rangle
(\Gamma_i^{\alpha}-\Gamma_{i-1}^{\alpha})
\right)=E\langle\varepsilon_i,x\rangle^2
\sum_{i=1}^n\E(\Gamma_i^{\alpha}-\Gamma_{i-1}^{\alpha})^2\sim C(x)n^{2\alpha-1},\,\, n\rightarrow\infty,
\]
where $C(x)=\frac{2\alpha^2}{2\alpha-1}E\langle\varepsilon_1,x\rangle^2.$
Therefore it is natural to take $B_n^2=n^{2\alpha-1}$.
\vspace{5pt}

We proceed in 5 steps:
\vspace{5pt}

{\bf Step 1}: Lemmas
\vspace{5pt}

{\bf Step 2}: We compare $X_n(\cdot)$ with $Z_n(\cdot)$ where $Z_n(t_{n,k})=\frac{\alpha}{B_n}\sum_{i=1}^k\varepsilon_i\gamma_i\Gamma_{i-1}^{\alpha-1}$ and show that $\|X_n-Z_n\|_{\infty}\xrightarrow{\mathbb{P}}0$.
\vspace{5pt}

{\bf Step 3}: We compare $Z_n(\cdot)$ with $W_n(\cdot)$ where $W_n(t_{n,k})=\frac{\alpha}{B_n}\sum_{i=1}^k\varepsilon_i\gamma_i(i-1)^{\alpha-1}$ and state that $\|Z_n-W_n\|_{\infty}\xrightarrow{\mathbb{P}}0$.
\vspace{5pt}

{\bf Step 4}: We show that process $U_n(\cdot)$,
\[
U_n\left(\left(\frac{k}{n}\right)^{\alpha}\right)=\frac{\alpha}{B_n}\sum_{i=1}^k\varepsilon_i\gamma_i(i-1)^{\alpha-1},
\]
converges weakly to the limiting process
\[
Y(t)=\sqrt{2\alpha}\int_0^ts^{\frac{\alpha-1}{2\alpha}}\,\mathrm{d}W(s);
\]
here $W(\cdot)$ is a  process of Brownian motion, for which the covariance matrix of $W(1)$
coincides with  the covariance matrix of $\varepsilon _{1}.$
\vspace{5pt}

{\bf Step 5}: We show that the convergence $W_n\Rightarrow Y$ follows from the convergence $U_n\Rightarrow Y$.
\vspace{5pt}

{\bf Finally}: We get the convergence $X_n\Rightarrow Y$.

\subsubsection{Step 1}
\begin{lemma}
Let $\alpha > 0$ and  $m \ge 1$. Then $\forall x>0,\,h>0$
\begin{equation}\label{eq:3}
(x+h)^{\alpha}-x^{\alpha}=\sum_{k=1}^ma_kh^kx^{\alpha-k}+R(x,h),
\end{equation}
where
\[
a_k=\frac{\alpha(\alpha-1)\dots(\alpha-k+1)}{k!},
\]
and
\begin{equation}
|R(x,h)|\le |a_{m+1}|h^{m+1} \max \{x^{\alpha-(m+1)}, (x+h)^{\alpha-(m+1)}\}.
\end{equation}
\end{lemma}
%\begin{proof}
{\bf Proof.}
By the formula of Taylor-Lagrange we have (\ref{eq:3}) with
\[
|R(x,y)|\le\frac{1}{(m+1)!}h^{m+1}\sup_{x\le t\le x+h}|f^{(m+1)}(t)|,
\]
where $f(t)=t^{\alpha}$. 
As $f^{(m+1)}(t)=\alpha(\alpha-1)\dots(\alpha-m)t^{\alpha-(m+1)}$, we get the result.

\halmos

\begin{lemma}
For $\alpha \geq 0$ and  $k\rightarrow \infty$
\begin{equation}
\left(1+\frac{\alpha}{k}\right)^k=e^\alpha+O\left(\frac{1}{k}\right).
\end{equation}
\end{lemma}
%\begin{proof}
{\bf Proof.}
It follows from the inequalities:
\[
0\le e^\alpha-\left(1+\frac{\alpha}{k}\right)^k\le\frac{e^{\alpha}\alpha^2}{k}.
\]
%\end{proof}
\halmos
\begin{lemma}
Let $\mathbf{\Gamma}$ be the Gamma function. Then as $k\rightarrow\infty$
\[
\frac{\mathbf{\Gamma}(k+\alpha)}{\mathbf{\Gamma}(k)}=k^{\alpha}+O(k^{\alpha-1}).
\]
\end{lemma}
%\begin{proof}
{\bf Proof.}
It follows from Lemma 2 and well known asymptotic
\[
\mathbf{\Gamma}(t)=t^{t-\frac{1}{2}}e^{-t}\sqrt{2\pi}\left(1+\frac{1}{12t}+O\left(\frac{1}{t^2}\right)\right),\,\,t\rightarrow\infty.
\]
%\end{proof}
\halmos

\begin{lemma}
For any real $\beta$
we have as $k\rightarrow\infty$
\[
E(\Gamma_k^{\beta})=k^{\beta}+O(k^{\beta-1}).
\]
\end{lemma}

{\bf Proof.}
The result follows from the well known fact that 
$$
E(\Gamma_k^{\beta})= 
\frac{\mathbf{\Gamma}(k+\beta)}{\mathbf{\Gamma}(k)}
$$
and Lemma 3.
\halmos

\begin{lemma} Let $\alpha \geq 0.$ 
The following relations take place as $k \rightarrow \infty$:

\begin{equation}\label{l5-1}
\Gamma_{k+1}^{\alpha}-\Gamma_k^{\alpha}=
\alpha\gamma_{k+1}\Gamma_k^{\alpha-1}+
\rho_k,
\end{equation}
where $|\rho_k| = O(k^{\alpha-2})$ in probability;
\vspace{5pt}

\begin{equation}\label{l5-2}
E|\Gamma_{k+1}^{\alpha}-\Gamma_k^{\alpha}|^2
= 2\alpha^2k^{2\alpha-2}+ O(k^{2\alpha-3})
;
\end{equation}

\begin{equation}\label{l5-3}
E|\Gamma_{k+1}^{\alpha}-\Gamma_k^{\alpha}-
\alpha\gamma_{k+1}\Gamma_k^{\alpha-1}|^2 =
O(k^{2\alpha-4}).
\end{equation}

\end{lemma}

We deduce immediately from (\ref{l5-2}) the following relation.
\begin{corollary}
We have
$$
\sum_1^{n-1}
E|\Gamma_{k+1}^{\alpha}-\Gamma_k^{\alpha}|^2
= \frac{2\alpha^2}{2\alpha-1}n^{2\alpha-1}+
O(n^{2\alpha-2}).
$$
\end{corollary}

{\bf Proof of Lemma 5.}
We find, applying Lemma 1,
\begin{equation}\label{L5}
\Gamma_{k+1}^{\alpha}-\Gamma_k^{\alpha}=
\alpha\gamma_{k+1}\Gamma_k^{\alpha-1}+
R(\Gamma_k,\gamma_{k+1}),
\end{equation}

where 
\begin{equation}\label{L6}
R(\Gamma_k,\gamma_{k+1}) \leq
\frac{1}{2}\gamma_{k+1}^2
\max_{\Gamma_k\leq s \leq\Gamma_{k+1}}
|\alpha(\alpha -1)| s^{\alpha-2}\leq
\frac{|\alpha(\alpha -1)|}{2}\gamma_{k+1}^2
\max\{\Gamma_{k+1}^{\alpha-2}, \Gamma_{k}^{\alpha-2}\} .
\end{equation}

As $\Gamma_{k} \sim k$ a.s. when $k\rightarrow \infty,$ we get (\ref{l5-1}).

The proofs of (\ref{l5-2}) and (\ref{l5-3})
follow directly from  (\ref{L5}), (\ref{L6}) and Lemma 4.
\halmos

\subsubsection{Step 2}
 We show that $\|X_n-Z_n\|_{\infty}\xrightarrow{\mathbb{P}}0$,
where
\[
Z_n(t_{n,k})=\frac{\alpha}{B_n}\sum_{i=1}^k\varepsilon_i\gamma_i\Gamma_{i-1}^{\alpha-1}.
\]

It is clear that
\[
\begin{aligned}
&\delta_n : =\|X_n-Z_n\|_{\infty}=\sup_{t\in[0,1]}|X_n(t)-Z_n(t)|=
\max_{k\le n}|X(t_{n,k})-Z_n(t_{n,k})|=\max_{k\le n}|r_k|,
\end{aligned}
\]
where 
$$
r_k=\frac{1}{B_n}\sum_{i=1}^k\varepsilon_i\left[\Gamma_i^{\alpha}-\Gamma_{i-1}^{\alpha}-
\alpha\gamma_i\Gamma_{i-1}^{\alpha-1}\right]=\sum_{i=1}^k\varepsilon_i\xi_i,
$$
and
  $$
  \xi_i=\left(\Gamma_i^{\alpha}-\Gamma_{i-1}^{\alpha}-\alpha\gamma_i\Gamma_{i-1}^{\alpha-1}\right)\frac{1}{B_n}.
  $$

Let $\mathfrak{M}=\sigma(\xi_1,\xi_2,\dots,\xi_n)=\sigma(\gamma_1,\gamma_2,\dots,\gamma_n).$ Under condition $\mathfrak{M}$ the sequence $(r_k)$ is the sequence of sums of independent random variables with mean zero. By Kolmogorov's inequality
\begin{equation}\label{eq:6}
\begin{aligned}
&\mathbb{P}\{\max_{k\le n}|r_k|\ge t\}=\E\{\mathbb{P}\{\max_{k\le n}|r_k|\ge t\,|\,\mathfrak{M}\}\}\le
\E\left(\frac{1}{t^2}\sum_{j=1}^n\xi^2_j\right)=\frac{1}{t^2}\sum_{j=1}^n\E\xi_j^2.
\end{aligned}
\end{equation}

By Lemma 5 $\E\xi_j^2 = O(j^{-3}).$
Therefore, 
\[
\sum_{j=1}^n\E\xi_j^2
= O(n^ {-2})
\]

Finally we get from (\ref{eq:6}): $\forall t>0$
\[
\mathbb{P}\{\delta_n\ge t\}\xrightarrow[n\rightarrow\infty]{}0,
\]
which gives the convergence $\|X_n-Z_n\|_{\infty}\xrightarrow{\mathbb{P}}0$.
\vspace{5pt}

\subsubsection{ Step 3} 
We show now that $\|Z_n-W_n\|_{\infty}\xrightarrow[n\rightarrow\infty]{\mathbb{P}}0$; where $W_n(t_{n,k})=\frac{\alpha}{B_n}\sum_{i=1}^k\varepsilon_i\gamma_i(i-1)^{\alpha-1}$.

We have
\[
\Delta_n=\sup_{t\in[0,1]}|Z_n(t)-W_n(t)|=\max_{k\le n}|Z_n(t_{n,k})-W_n(t_{n,k})|=\max_{k\le n}\{|\beta_k|\},
\]
where $\beta_k=\frac{\alpha}{B_n}\sum_{i=1}^k\varepsilon_i\gamma_i\left(\Gamma_{i-1}^{\alpha-1}-(i-1)^{\alpha-1}\right)$.

Similar to the previous case $(\beta_k)$ under condition $\mathfrak{M}$ is the sequence of sums of independent random variables with mean zero. Therefore
\[
\begin{aligned}
\mathbb{P}\{\max_{k\le n}\{|\beta_k|\}\ge t\}=\E\left(\mathbb{P}\{\max_{k\le n}\{|\beta_k|\}\ge t\,|\,\mathfrak{M}\}\right)\le\frac{1}{t^2}\sum_{j=1}^n\E\eta_j^2,
\end{aligned}
\]
where $\eta_j=\frac{\alpha}{B_n}\gamma_j\left(\Gamma_{j-1}^{\alpha-1}-(j-1)^{\alpha-1}\right)$.
\vspace{5pt}

{\bf Estimation of $\E\eta_j^2.$}

By independence of $\gamma_j$ and $\Gamma_{j-1}$
\[
\E\eta_j^2=\frac{2\alpha^2}{B_n^2}\E\left(\Gamma_{j-1}^{\alpha-1}-(j-1)^{\alpha-1}\right)^2
\]
Let us change $j-1$ to $k$
\[
\begin{aligned}
&\E\left(\Gamma_{k}^{\alpha-1}-k^{\alpha-1}\right)^2=\E\left(\Gamma_k^{2\alpha-2}\right)+k^{2\alpha-2}-2k^{\alpha-1}\E\left(\Gamma_k^{\alpha-1}\right)=\\
&=\frac{\Gamma(k+2\alpha-2)}{\Gamma(k)}+k^{2\alpha-2}-2k^{\alpha-1}\frac{\Gamma(k+\alpha-1)}{\Gamma(k)}=(\text{by Lemma\,3})=\\
&=\left[k^{2\alpha-2}+O(k^{2\alpha-3})+k^{2\alpha-2}-2k^{2\alpha-2}\right]=O(k^{2\alpha-3}).
\end{aligned}
\]

Hence
\[
\E\eta_j^2\le C\frac{j^{2\alpha-3}}{n^{2\alpha-1}}
\]
and
\[
\sum_{j=1}^n\E\eta_j^2\le C\frac{1}{n^2}
\]

We have finally $\mathbb{P}\{\max_{k\le n}|\beta_k|\ge t\}\rightarrow0,\,\,n\rightarrow\infty,$ which gives the convergence\\ $\|W_n-Z_n\|\xrightarrow{\mathbb{P}}0$.
\vspace{7pt}

\subsubsection{ Step 4} 
Let $U_n$ be the process defined at the points $\frac{k}{n}$ by
\[
U_n\left(\left(\frac{k}{n}\right)^{\alpha}\right)=\frac{\alpha}{B_n}\sum_{i=1}^k\varepsilon_i\gamma_i(i-1)^{\alpha-1},\,\,k=1,2,\dots,n,
\]
and by linear interpolation on the intervals $[\frac{k}{n}, \frac{k+1}{n}],\; k=0,\ldots,n-1.$
We now state weak convergence of the processes $U_n$ 
to the process $Y$,
\[
Y(t)=\sqrt{2\alpha}\int_0^ts^{\frac{\alpha-1}{2\alpha}}\,\mathrm{d}W(s),
\]
W is a Brownian motion, for which the covariance matrix of 
$W(1)$
coincides with  the covariance matrix of $\varepsilon _{1}.$

The proof is standard because $U_n(\cdot)$ represents a (more or less) usual broken line constructed by the consecutive sums of independent (non-identically distributed) random variables. One could apply Prokhorov's theorem
(see Gikhman and Skorohod (1996), ch.IX, sec. 3, Th.1).

Only one thing must be checked: that for any $0<s<t\le1,\;\;$
and for any $x\in {\mathbb R}^d,\; |x|=1,$ we have the convergence
 $\langle U_n(t)-U_n(s), x \rangle\Longrightarrow \langle Y(t)-Y(s), x \rangle$.

It is clear that
\[
[U_n(t)-U_n(s)]-\left[U_n\left(\left(\frac{k}{n}\right)^{\alpha}\right)-U_n\left(\left(\frac{l}{n}\right)^{\alpha}\right)\right]\;\;\xrightarrow{\;\;\mathbb{P}\;\;}\;\;0,
\]
if $\left(\frac{k}{n}\right)^{\alpha}\rightarrow t,\;\;\; \left(\frac{l}{n}\right)^{\alpha}\rightarrow s$.
\vspace{7pt}

Let $l <k.$ As 
$$
\left\langle U_n\left(\left(\frac{k}{n}\right)^{\alpha}\right)-U_n\left(\left(\frac{l}{n}\right)^{\alpha}\right), x 
\right\rangle=\frac{\alpha}{B_n}\sum_{i=l+1}^k\langle\varepsilon_i, x \rangle\gamma_i(i-1)^{\alpha-1},
$$
by the theorem of Lindeberg-Feller it is sufficient to state the convergence of variances.

We have
\[
\Var\left\langle U_n\left(\left(\frac{k}{n}\right)^{\alpha}\right)-U_n\left(\left(\frac{l}{n}\right)^{\alpha}\right), x 
\right\rangle
=
\]
\[
= \frac{2\alpha^2}{n^{2\alpha-1}}
E\langle\varepsilon_1, x \rangle^2\sum_{i=l+1}^k(i-1)^{2\alpha-2}\xrightarrow[n\rightarrow\infty]{}\frac{2\alpha^2}{2\alpha-1}E\langle\varepsilon_1, x \rangle^2[t^{\frac{2\alpha-1}{\alpha}}-s^{\frac{2\alpha-1}{\alpha}}],
\]
and
\[
\Var\langle Y(t)-Y(s), x\rangle=2\alpha E\langle\varepsilon_1, x \rangle^2\int_s^tu^{\frac{\alpha-1}{\alpha}}\,\mathrm{d}u=\frac{2\alpha^2}{2\alpha-1}E\langle\varepsilon_1, x \rangle^2[t^{\frac{2\alpha-1}{\alpha}}-s^{\frac{2\alpha-1}{\alpha}}],
\]
which are the same.

\subsubsection{ Step 5:  Convergence $X_n\Rightarrow Y$.}

Due to the steps 2 and 3 it is sufficient to show that $W_n\Rightarrow Y$.

Let $f_n:\,[0,1]\rightarrow[0,1]$, be a piecewise linear continuous function such that $f_n(t_{n,k})=\left(\frac{k}{n}\right)^{\alpha}$; $t_{n,k}=\left(\frac{\Gamma_k}{\Gamma_n}\right)^{\alpha}$; $k=0,1,\dots,n$.

By definition of $W_n$ and $U_n$ we have
\[
W_n(t)=U_n(f_n(t)), \;t \in [0,1].
\]

By the corollary to Lemma \ref{lem:6} (see below) the function $f_n$  converges in probability uniformly to  $f$, $f(t)=t$, and by previous step $U_n \Rightarrow Y$.

It means that we can apply Lemma \ref{lem:7} which gives the necessary convergence.
\begin{lemma}\label{lem:6}
Let 
$$
M_n=\max_{k\le n}\left\{ \left|\frac{\Gamma_k}{\Gamma_n}-\frac{k}{n}\right|\right\}.
$$ 

Then $\;\;M_n\xrightarrow{\mathbb{P}}0,\;\;\;n\rightarrow\infty$.
\end{lemma}
{\bf Proof of Lemma 6.}
We have
\begin{equation}\label{eq:7}
\begin{aligned}
&\mathbb{P}\{M_n>\varepsilon\}=\E\left\{\mathbb{P}\left\{\max_{k\le n}\left|\frac{\Gamma_k}{\Gamma_n}-\frac{k}{n}\right|>\varepsilon\, |\,\Gamma_n \right\}\right\}=\\
&=\int_{0}^{\infty}\mathbb{P}\left\{\max_{k\le n}\left|\frac{\Gamma_k}{\Gamma_n}-\frac{k}{n}\right|>\varepsilon\, |\,\Gamma_n=t \right\}\mathcal{P}_{\Gamma_n}(\mathrm{d}t)=\\
&=\int_{0}^{\infty}\mathbb{P}\left\{\max_{k\le n}\left|\xi_{n,k}-\frac{k}{n}\right|>\varepsilon \right\}\mathcal{P}_{\Gamma_n}(\mathrm{d}t)
= \mathbb{P}\left\{\max_{k\le n}\left|\xi_{n,k}-\frac{k}{n}\right|>\varepsilon \right\},
\end{aligned}
\end{equation}
where $(\xi_{n,k})_{k=1,\dots,n}$ are the order statistics from $[0,1]  $-uniform distribution.

Let $\delta_n : =\max_{k\le n}\left|\xi_{n,k}-\frac{k}{n}\right|$. Evidently, $\delta_n\le \sup_{[0,1]}|F_n^*(x)-x|$, where $F_n^*$ is a uniform empirical distribution function. By Glivenko-Cantelli theorem, $\sup_{[0,1]}|F_n^*(x)-x|\rightarrow0$ a.s, which gives the convergence $M_n \rightarrow0$ in probability.
%\end{proof}
\halmos
\begin{corollary}
$M_{n}^{(1)}=\max_{k\le n}\left|\left(\frac{\Gamma_k}{\Gamma_n}\right)^{\alpha}-\left(\frac{k}{n}\right)^{\alpha}\right|\xrightarrow{\mathbb{P}}0$, $n\rightarrow\infty$.
\end{corollary}

The proof follows directly from Lemma \ref{lem:6} due to the uniform continuity of the function $h(x)=x^{\alpha}$, $x\in[0,1]$.

\begin{lemma}\label{lem:7}
Let $\{U_n\}$ be a sequence of continuous processes on $[0,1]$ weakly convergent to some limit process $U$. Let $\{f_n\}$ be a sequence of random continuous bijections $[0,1]$ on $[0,1]$ which in probability uniformly converges to the identity function $f(t)\equiv t$. Then the process $W_n,\,\, W_n(t)=U_n(f_n(t)),\, t\in[0,1],$ will converge weakly to $U$.
\end{lemma}

{\bf Proof of Lemma 7.}
By theorem 4.4 from Billingsley (1968) we have the weak convergence in 
${\mathbb M} : ={\mathbb C}[0,1] \times {\mathbb C}[0,1]$
$$
(U_n,f_n) \Longrightarrow (U,f).
$$
By Skorohod representation theorem we can find a random elements $(\tilde U_n,\tilde f_n)$ and $  (\tilde U,\tilde f)$
of ${\mathbb M}$ (defined probably on a new probability space) such that
$$
(U_n,f_n) \stackrel{\cal L}{=}(\tilde U_n,\tilde f_n),\;\;\;(U,f) \stackrel{\cal L}{=} (\tilde U,\tilde f),
$$
and
$(\tilde U_n,\tilde f_n) \rightarrow  (\tilde U,\tilde f)$ a.s. in ${\mathbb M}.$

As the last convergence implies evidently the a.s. uniform convergence of 
$\tilde U_n(\tilde f_n(t))$ to $\tilde U(\tilde f(t)),$
we get the convergence in distribution of $U(f_n(\cdot))$ to $ U(f(\cdot)) = U(\cdot).$
\halmos

\section{Proof of Th. 2}
{\bf Proof of Th. 2.} We need some facts from  Stroock and  Varadhan
(1979). Consider $%
(\Omega ,\mathcal{M}),$ where $\Omega ={\mathbb C}([0,\infty );R^{d})$ be the space of
continuous trajectories from $[0,\infty )$ into $R^{d}.$ Given $t\geq 0$ and $%
\omega \in \Omega $ let $x(t,\omega )$ denote the position of $\omega $ in $%
R^{d}$ at time $t.$ If we put%
\[
D(\omega ,\omega ^{\prime })=\sum_{n=1}^{\infty }\frac{1}{2^{n}}\frac{%
\sup_{0\leq t\leq n}\left\vert x(t,\omega )-x(t,\omega ^{\prime
})\right\vert }{1+\sup_{0\leq t\leq n}\left\vert x(t,\omega )-x(t,\omega
^{\prime })\right\vert }
\]%
then it is well known that $D$ is a metric on $\Omega $ and $\left( \Omega
,D\right) $ is a Polish space. The convergence induced by $D$ is uniform
convergence on bounded $t$ - intervals. 
For simplicity, we will omit $\omega$ in the future and we will
 be assuming that all
our processes are homogeneous in time. Analogous results for
time-inhomogeneous processes may be obtained by simply considering the
time-space processes. 

We will use $\mathcal{M}$ to denote
the Borel $\sigma $ - field of subsets of  $\left( \Omega ,D\right) ,%
\;\mathcal{M=\sigma \lbrack }x(t):t\geq 0].$ We also will consider an
increasing family of $\sigma $-algebras  $\mathcal{M}_{t}=\mathcal{\sigma
\lbrack }x(s):0\leq s\leq t].$ 
Classical approach to the construction of diffusion
processes corresponding to given coefficients $a$ and $b$ involves a
transition probability function $P(s,x;t,\cdot )$ which allows to construct 
for each $x$ $\in R^{d},$ a probability measure $P_{x}$ on $\Omega
={\mathbb C}([0,\infty );R^{d})$ with the properties that%
\[
P_{x}(x(0)=x)=1
\]%
and%
\[
P_{x}(x(t_{2})\in \Gamma \text{ }|\mathcal{M}%
_{t_{1}})=P(t_{1},x(t_{1});t_{2},\Gamma )\text{ }a.s.P_{x}
\]%
for all $0\leq t_{1}<t_{2}$ and $\Gamma \in \mathcal{B}_{R^{d}}$ (the Borel $%
\sigma $ - algebra in $R^{d}).$ It appears that this measure is a martingale
measure for a special martingale related with the second order differential
operator 
\[
L=\frac{1}{2}\sum_{i,j=1}^{d}a^{ij}(\cdot )\frac{\partial ^{2}}{\partial
x_{i}\partial x_{j}}+\sum_{i=1}^{d}b^{i}(\cdot )\frac{\partial }{\partial
x_{i}},
\]%
namely, for all $f\in {\mathbb C}_{0}^{\infty }(R^{d})$%
\[
P_{x}(x(0)=x)=1,
\]
\begin{equation}
(f(x(t))-\int_{0}^{t}Lf(x(u))du,\mathcal{M}_{t},P_{x})  \label{d1}
\end{equation}%
\bigskip is a martingale. We will say that the martingale problem for $a$
and $b$ is \textit{well-posed \ }if, for each $x$ there is exactly one
solution to that martingale problem starting from $x.$ We will be working
with the following set up. For each $h>0$ let $\Pi _{h}(x,\cdot )$ be a
transition function on $R^{d}.$ Given $x\in R^{d},$ \,let $P_{x}^{h}$ be the
probability measure on $\Omega $ characterized by the properties that%
\begin{equation}
(i)\text{ \ \ }P_{x}^{h}(x(0)=x)=1,  \label{d2}
\end{equation}%
\begin{equation}\label{d3}
(ii)\text{ \ \ }P_{x}^{h}\Big\{ x(t)=\frac{(k+1)h-t}{h}x(kh)+\frac{t-kh}{h}%
x((k+1)h),\;\;  kh\leq t<(k+1)h\Big\} =1
\end{equation}
$$
\text{ for all }k\geq 0,  
$$

\[
(iii)\text{ \ }P_{x}^{h}(x((k+1)h)\in \Gamma \text{ }|\text{ }\mathcal{M}%
_{kh})=\Pi _{h}(x(kh),\Gamma ),\;\;\;P_{x}^{h}-\text{ a.s.}
\]%
\begin{equation}
\text{for all }k\geq 0\text{ and }\Gamma \in \mathcal{B}_{R^{d}}.  \label{d4}
\end{equation}

Define%
\begin{equation}
a_{h}^{ij}(x)=\frac{1}{h}\int_{\left\vert y-x\right\vert \leq
1}(y_{i}-x_{i})(y_{j}-x_{j})\Pi _{h}(x,dy),  \label{d5}
\end{equation}%
\begin{equation}
b_{h}^{i}(x)=\frac{1}{h}\int_{\left\vert y-x\right\vert \leq
1}(y_{i}-x_{i})\Pi _{h}(x,dy),  \label{d6}
\end{equation}%
and%
\begin{equation}
\Delta _{h}^{\varepsilon }(x)=\frac{1}{h}\Pi _{h}(x,R^{d}\backslash
B(x,\varepsilon )),  \label{d7}
\end{equation}%
where $B(x,\varepsilon )$ is the open ball with center $x$ and radius $\varepsilon .$
What we are going to assume is that for all $R>0$%
\begin{equation}
\lim_{h\searrow 0}\sup_{\left\vert x\right\vert \leq R}\left\Vert
a_{h}(x)-a(x)\right\Vert =0,  \label{d8}
\end{equation}%
\begin{equation}
\lim_{h\searrow 0}\sup_{\left\vert x\right\vert \leq R}\left\vert
b_{h}(x)-b(x)\right\vert =0,  \label{d9}
\end{equation}%
\begin{equation}
\sup_{h>0}\sup_{x\in R^{d}}\left( \left\Vert a_{h}(x)\right\Vert +\left\vert
b_{h}(x)\right\vert \right) <\infty ,  \label{d10}
\end{equation}%
\begin{equation}
\lim_{h\searrow 0}\sup_{x\in R^{d}}\Delta _{h}^{\varepsilon }(x)=0.\text{ }
\label{d11}
\end{equation}%

\textbf{Theorem} \textbf{A}. (Strook and Varadhan (1979), page 272, Theorem 11.2.3). \textit{%
Assume that in addition to (\ref{d8})-(\ref{d11}) the coefficients }$a$ 
\textit{and} $b$ \textit{are continuous and have the property that for each }%
$x\in R^{d}$ \textit{the martingale problem for }$a$ \textit{and} $b$ 
\textit{has exactly one solution }$P_{x}$ \textit{starting from }$x$ \textit{%
(that is well posed). Then  }$P_{x}^{h}$ converges weakly to  $P_{x}$ \textit{%
uniformly in }$x$ \textit{on compact subsets of }$R^{d}.$

Sufficient conditions for the well-posedness
%incorrect !!!
is given by the following
theorem.
\vspace{7pt}

Let $S_d$ be the set of symmetric non-negative definite $d\times d$ real matrices.
\vspace{7pt}

\textbf{Theorem B. }(Strook and Varadhan (1979), page 152, Theorem 6.3.4). \textit{Let }$%
a:R^{d}\longrightarrow S_{d}$ \textit{and }$b:R^{d}\longrightarrow R^{d}$ \ 
\textit{be bounded measurable functions and suppose that }$\sigma
:R^{d}\longrightarrow R^{d}\times R^{d}$ \textit{is a bounded measurable
function \ such that \ }$a=\sigma \sigma ^{\ast }.$ \textit{Assume that
there is an }$A$ \textit{such that}
\vspace{10pt}

\begin{equation}
\left\Vert \sigma (x)-\sigma (y)\right\Vert +\left\vert b(x)-b(y)\right\vert
\leq A\left\vert x-y\right\vert   \label{d12}
\end{equation}%
\textit{for all \ }$x,y\in R^{d}.$ \textit{Then the martingale problem for }$%
a$ \textit{and} $b$ \textit{is well-posed and the corresponding family of
solutions }$\{P_{x}:x\in R^{d}\}$ \textit{is Feller continuous (that is }$%
P_{x_{n}}\rightarrow P_{x}$ \textit{weakly if }$x_{n}\rightarrow x$\textit{).%
}

Note that (\ref{d12}) and uniform ellipticity of $a(x)$ imply the existence
of the transition density $p(s,x;t,y)$ (Strook and Varadhan (1979), Theorem 3.2.1, page 71). 

Consider the model%
\[
X\left( (k+1)h\right) =X(kh)+hb(X(kh))+\sqrt{h}\xi (X(kh)),
\]%
\begin{equation}
\xi (X(kh))=\rho _{k}\sigma (X(kh))\varepsilon _{k},  \label{d13}
\end{equation}%
where $\left\{ \varepsilon _{k}\right\} $ are i.i.d. random vectors uniformly
distributed on the unit sphere $S^{d-1},$ and $\left\{ \rho _{k}\right\} $
are i.i.d. random variables having a density, $\rho _{k}\geq 0,$ $E\rho
_{k}^{2}=d$. Let us check the conditions (\ref{d8})-(\ref{d11}). It is easy to
see that \ 
\begin{equation}
\Pi _{h}(x,dy)=p_{h}^{x}(y)dy,\text{ where }\;p_{h}^{x}(y)=h^{-d/2}f_{\xi }\left(%
\frac{y-x-hb(x)}{\sqrt{h}}\right). \label{d14}
\end{equation}%
Here $f_{\xi }$ denotes the density of the random vector $\xi .$ Let us check
(\ref{d11}). Note that $E\xi =0$ and the covariance matrix of the vector $\xi 
$ is equal to 
\begin{equation}
Cov(\xi ,\xi ^{T})=E(\rho _{k}^{2}\sigma (x)\varepsilon _{k}\varepsilon
_{k}^{T}\sigma ^{T}(x))=a(x).  \label{d15}
\end{equation}%
We have%
\[
h\Delta _{h}^{\varepsilon }(x)=\Pi _{h}(x,R^{d}\backslash B(x,\varepsilon
))=\int_{R^{d}\backslash B(x,\varepsilon )}p_{h}^{x}(y)dy=
\]%
\[
=\int_{v+\sqrt{h}b(x)\in R^{d}\backslash B(0,\frac{\varepsilon }{\sqrt{h}}%
)}f_{\xi }(v)dv=P\left\{\xi \in \overline{B\left(0,\frac{\varepsilon }{\sqrt{h}}\right)}\right\}-%
\sqrt{h}b(x))\;\leq 
\]%
\begin{equation}
\leq \;P\left\{\left\vert \xi \right\vert ^{2}\geq \frac{\varepsilon ^{2}}{4h}\right\}=o(h).
\label{d16}
\end{equation}%
The last equality is a consequence of the Markov inequality. The equality (\ref{d15}),
uniform ellipticity of $a(x)$ and (\ref{d16}) implies (\ref{d11}). To prove (\ref%
{d8}) note that by (\ref{d12})%
\[
a_{h}^{ij}(x)=\frac{1}{h}\int_{\left\vert y-x\right\vert \leq
1}(y_{i}-x_{i})(y_{j}-x_{j})p_{h}^{x}(y)dy=
\]%
\[
=\int_{\left\vert v+\sqrt{h}b(x)\right\vert \leq \frac{1}{\sqrt{h}}}(v_{i}+%
\sqrt{h}b^{i}(x))(v_{j}+\sqrt{h}b^{j}(x))f_{\xi }(v)dv=
\]%
\begin{equation}
=\int_{\left\vert v+\sqrt{h}b(x)\right\vert \leq \frac{1}{\sqrt{h}}%
}v_{i}v_{j}f_{\xi }(v)dv+o(\sqrt{h})=a(x)+o(1).  \label{d17}
\end{equation}%
To check (\ref{d9}) note that 
\[
b_{h}^{i}(x)=\frac{1}{h}\int_{\left\vert y-x\right\vert \leq
1}(y_{i}-x_{i})p_{h}^{x}(y)dy=
\]%
\[
=\frac{1}{\sqrt{h}}\int_{\left\vert v+\sqrt{h}b(x)\right\vert \leq \frac{1}{%
\sqrt{h}}}(v_{i}+\sqrt{h}b^{i}(x))f_{\xi }(v)dv=
\]%
\begin{equation}
=b^{i}(x)\int_{\left\vert v+\sqrt{h}b(x)\right\vert \leq \frac{1}{\sqrt{h}}%
}f_{\xi }(v)dv-\frac{1}{\sqrt{h}}\int_{\left\vert v+\sqrt{h}b(x)\right\vert >%
\frac{1}{\sqrt{h}}}v_{i}f_{\xi }(v)dv.  \label{d18}
\end{equation}%
To estimate the second integral in (\ref{d18}) we apply the Cauchy - Schwarz
inequality 
\begin{equation}
\frac{1}{\sqrt{h}}\int_{\left\vert v+\sqrt{h}b(x)\right\vert >\frac{1}{\sqrt{%
h}}}\left\vert v\right\vert f_{\xi }(v)dv\leq \frac{1}{\sqrt{h}}\left( \int
\left\vert v\right\vert ^{2}f_{\xi }(v)dv\right) ^{1/2}\left( P(\left\vert
\xi \right\vert ^{2}\geq \frac{1}{4h}\right) ^{1/2}=o(1),  \label{d19}
\end{equation}%
and (\ref{d18}), (\ref{d19}) \ imply (\ref{d9}). Finally,  (\ref{d10}) follows
from our calculations and assumptions of Theorem B. Weak convergence \textit{%
\ }$P_{x}^{h}$  to $P_{x}$ follows now from Theorems A and B cited above.

\halmos

%\section{Conclusion}

\section{References}

\bibliographystyle{plain}

\begin{enumerate}

\bibitem  {[B]} P. Billingsley, Convergence of probability measures,
1968, John Wiley and Sons, N-Y.

\bibitem  {[GS]} I. I. Gikhman and A. V. Skorohod, Introduction to the theory of random processes,
 1996,  Dover Publications. 

\bibitem  {[GR]} I. Gradshtein and I. Ryzhik. Tablicy integralov, summ, ryadov i proizvedenii (rus.). Fizmatgiz,
M., (1963)

\bibitem  {[K]} A. D. Kolesnik, The explicit probability distribution of
a six-dimensional random flight, Th. of Stoch. Proc., v.15 (30), 1 (2009),
pp. 33-39.

\bibitem  {[Ko]} V. Konakov. Metod parametriksa diya diffusii i cepei
Markova.   Preprint (in Russian). Izdatel'stvo popechitel'skogo soveta
mehaniko-matematiceskogo fakul'teta MGU. Seriya WP BRP "STI". 2012. 

\bibitem  {[Kl]} J. C. Kluyver, A local probability problem, in: Proceedings of the Section of Sciences, Koninklijke Akademie van
Wetenschappen te Amsterdam, vol. 8, 1905, pp. 341–350.

\bibitem  {[KoM]} V. Konakov, E. Mammen. Small time Edgeworth-type expansions for weakly
convergent non homogenious Markov chains, PTRF, 143, (2009), 137-176.

\bibitem{[M]}  B. Mandelbrot, The fractional geometry of Nature, N-Y
(1982).

\bibitem {[O1]} E. Orsingher, A. De Gregorio, Reflecting Random Flights,\\
J. of Stat.Phys. \ 160 (6), (2015), pp. 1483-1506.

\bibitem {[O2]} E. Orsingher, R. Garra, Random flights governed by
Klein-Gordon type partial differential equations, Stoch. Proc. and
Appl., v.124 (2014), pp. 2171-2187.

\bibitem {[O3]} E. Orsingher, A. De Grigorio, Flying randomly in $R^{d}$
with Dirichlet displacements, Stoch. Proc. and Appl., (2012), v.122, pp. 676-713.

\bibitem  {[R]} L. Rayleigh, On the problem of the random flights and of random vibrations in one, two and three dimensions,
Philosophical Magazine, 37 (1919), 321–347.


\bibitem  {[SV]} D. W. Stroock and S. R. S. Varadhan, Multidimensional Diffusion Processes,
1979, Springer-Verlag.

\bibitem  {[Ya]} M. I. Yadrenko, Spektral'naya teoriya sluchainyh
polei (rus.), Kiev, (1980).

\end{enumerate}

\end{document}